\documentclass{nature}
\usepackage[margin=1in]{geometry}
\usepackage{mathpazo}
\usepackage{hyperref}
\hypersetup
{
	colorlinks,	%
	citecolor=black,%
	filecolor=black,%
	linkcolor=blue,%
	urlcolor=blue
}

\usepackage{amssymb,amsmath,epsfig,makecell,listings, amsthm, mathtools, paralist, enumitem, mathrsfs}
\usepackage{lineno}
\linenumbers
\usepackage{todonotes}
\newcounter{todocounter}

\presetkeys{todonotes}{inline, color=yellow!30}{}
\usepackage{nameref}
\usepackage[font=small, labelfont=bf]{caption}
\usepackage[font=footnotesize, labelfont=bf]{subcaption}

\usepackage{algorithm}
\usepackage{algorithmic}
\floatname{algorithm}{Procedure}

\newcommand{\R}{\ensuremath{\mathbb{R}}}
\newcommand{\Chn}{\ensuremath{\textbf{C}}}
\newcommand{\Hom}{\ensuremath{\textbf{H}}}
\newcommand{\VR}{\ensuremath{\underline{\mathrm{VR}}}}
\newcommand{\Barc}{\ensuremath{\mathrm{Bar}}}

\nolinenumbers

\title{Geometric anomaly detection in data} 
\author{Bernadette J Stolz$^{1}$, Jared Tanner$^{1,2}$, Heather A Harrington$^{1,2}$  \& Vidit Nanda$^{1,2}$}

\usepackage{graphicx}
\makeatletter
\let\saved@includegraphics\includegraphics
\AtBeginDocument{\let\includegraphics\saved@includegraphics}
\renewenvironment*{figure}{\@float{figure}}{\end@float}
\makeatother

\begin{document}

\maketitle

\begin{affiliations}
\item Mathematical Institute, University of Oxford, Oxford, UK
\item The Alan Turing Institute, London, UK
\end{affiliations}

\begin{abstract}
	
% word count: 149/200 [real max 300]
\medskip 

{This paper describes the systematic application of local topological methods for detecting interfaces and related anomalies in complicated high-dimensional data. By examining the topology of small regions around each point, one can optimally stratify a given dataset into clusters, each of which is in turn well-approximable by a suitable submanifold of the ambient space. Since these approximating submanifolds might have different dimensions, we are able to detect non-manifold like singular regions in data even when none of the data points have been sampled from those singularities. We showcase this method by identifying the intersection of two surfaces in the 24-dimensional space of cyclo-octane conformations, and by locating all the self-intersections of a Henneberg minimal surface immersed in 3-dimensional space. Due to the local nature of the required topological computations, the algorithmic burden of performing such data stratification is readily distributable across several processors.}
\end{abstract}

\begin{figure}[ht!]
	\centering
	{\includegraphics[width=.8\textwidth]{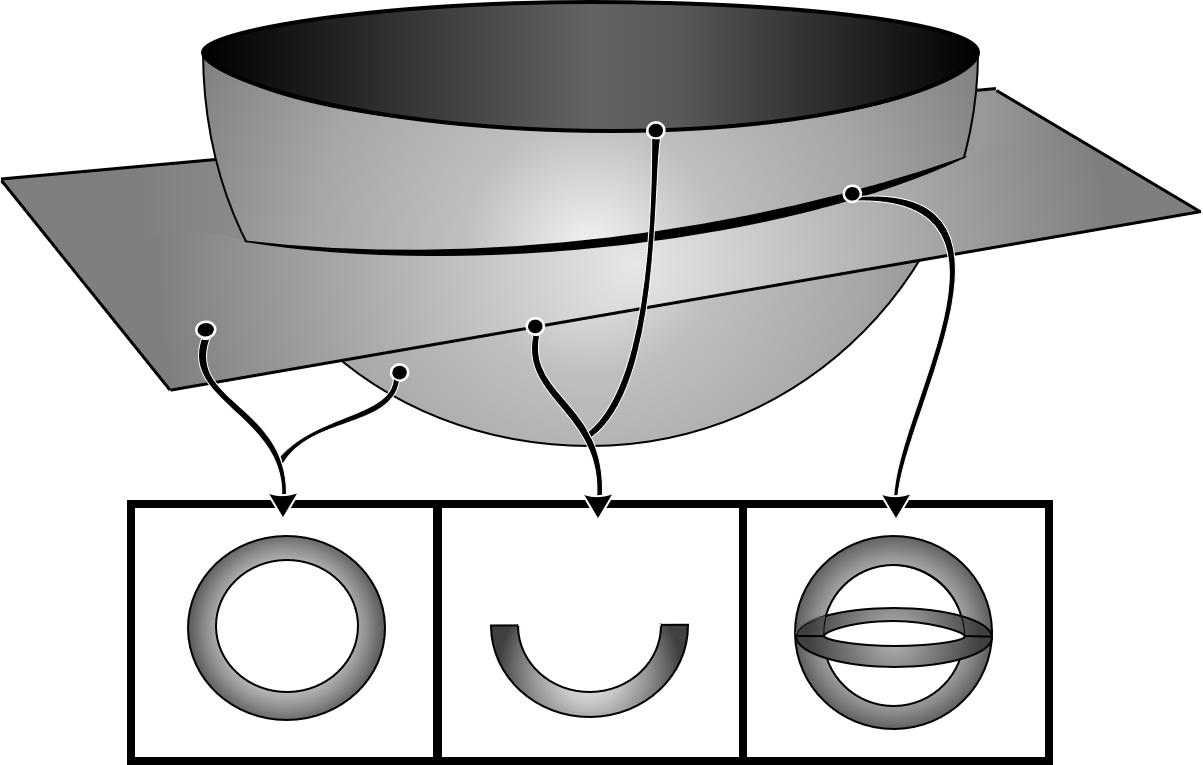}} %106 words
	\caption{Annular neighbourhood classes $A_x$ of several points $x$ in union of a hemisphere with a plane along an ellipse. All points lying far from this ellipse and from the boundaries have $A_x$ which look like a thickened circle, as shown in the left-most panel. All points lying in the boundary have $A_x$ which resembles a thickened half-circle, depicted in the middle panel. And all points $x$ on the singular ellipse itself have neighbourhoods $A_x$ which resemble two thickened circles glued along two edges, as in the rightmost panel. The dimensions of $\Hom^1(A_x)$ from left to right are 1, 0 and 3 respectively.}\label{Fig:Hemisphere_Plane}
\end{figure}

% 227 Words
The {\it manifold hypothesis}\cite{MR3522608} forms a cornerstone of modern data science; it asserts that the points in a naturally-occurring dataset tend to cluster near a manifold of dimension substantially lower than the ambient dimension of the data.  Typical examples of inferential schemes which rely on the manifold hypothesis include ({\bf a})
classical {\em principal component analysis}\cite{whatispca}, where data is approximated by an affine subspace, ({\bf b}) {\em visual perception}\cite{Seung2268}, where continuous changes of pose of an object yield smoothly-varying changes along a curved manifold, ({\bf c}) {\em subspace clustering}\cite{5714408}, where data is clustered into disjoint sets that are well approximated by affine subspaces, and ({\bf d}) {\em generative adversarial networks}, which naturally produce data on pairs of manifolds\cite{DBLP:conf/iclr/CheLJBL17}. In sharp contrast to this profusion, one encounters a remarkable dearth of techniques designed for the analysis of data sampled from non-manifold, or singular, spaces. Among the simplest examples of singular spaces are unions of two manifolds along a common submanifold (as shown in Figure \ref{Fig:Hemisphere_Plane}); these arise organically when more than one class of data are present in the same set of observations. Recent techniques for the analysis of such heterogeneous data, see for instance {\em capsule networks}\cite{NIPS2017_6975}, have focused primarily on coherently fusing together the multiple data classes.

% 177 words
The present work is motivated by an antipodal philosophy --- we believe that singular regions of spaces which underlie modern datasets are inherently interesting, that they will play an increasingly important role in the future of data analysis, and that it is therefore of paramount importance to be able to detect these singularities directly from the data points. Here we describe a new algorithm to accomplish this task --- in particular, we use a geometric approach to identify which data points lie near the intersection of more than one manifold. Our algorithm is based on {\em local cohomology}\cite{nanda2017local} and the theory of {\em stratifications}\cite{kirwan2006introduction}, which form particularly rich and fruitful enterprises in the study of singular spaces that arise in algebraic topology\cite{IH2} and geometry\cite{fulton}. Recent computational advances in these fields\cite{mischaikow:nanda:13,henselman:ghrist:16} have made it possible to bring this formidable theory to bear on the very concrete task of analysing data which lives on, or even near, spaces that are far more complicated than manifolds. 

% 162 words
Manifolds of dimension $n$ are characterized by the requirement that a small neighbourhood around each point should resemble the $n$-dimensional Euclidean disk (up to a standard equivalence relation called homeomorphism). While there can be no algorithmic procedure to determine whether two $n$-manifolds are homeomorphic or not\cite{markov} for $n > 4$, algebraic topology offers recourse to several rigorous descriptors for testing weaker forms of equivalence. Among the best known computable homeomorphism-invariants is {\em cohomology}, which assigns a sequence $\Hom^i(X)$ of vector spaces to a given topological space $X$. Although cohomology does not distinguish between Euclidean disks of different dimensions (all of these have the same cohomology as that of a point), it is an excellent tool for distinguishing $n$-dimensional spheres $\mathbb{S}^n$ from each other across different choices of $n$. Indeed, for all $n > 0$, we have 
\[
\dim \Hom^i(\mathbb{S}^n) = \begin{cases}
1 & \text{ if }i = n, \\
0 & \text{ otherwise.}
\end{cases}
\]
    
%153 words
Since the boundary of an $n$-dimensional disk is an $(n-1)$-dimensional sphere, our strategy for detecting singular regions in a dataset $P$ of points in Euclidean space $\R^n$ is as follows: we fix two real parameters $0 < r < s$, and around each point $x$ of $P$ we examine the subset of {\em annular neighbours} $A_x$ of $x$ --- this set consists of those points $y$ in $P$ whose Euclidean distance to $x$ satisfies $r \leq \|x-y\| \leq s$, and it forms a discrete proxy for the boundary of a neighbourhood around $x$. We then compute the cohomology of $A_x$ at various scales (often called the {\em persistent cohomology} of $A_x$), and use this information to quantify whether or not $A_x$ approximates a single sphere of some fixed dimension. If the answer is negative, then -- provided we have made judicious choices of $r$ and $s$ -- the point $x$ lies near a singular region of $X$.
Points which are singular due to the intersection of low-dimensional manifolds can alternatively be identified by measuring changes in the local dimension; for instance, the relative sum-of-squares of the for first $n-1$ singular values of $A_x$ is nearly one for $x$ away from an intersection\cite{Martin2011}.
%and can be as small as $1/2$ for $x$ near two orthogonally intersecting $n-1$ dimensional hyperplanes.

\begin{figure}[ht!]
	\centering
	\subcaptionbox{}{\includegraphics[width=.4\textwidth]{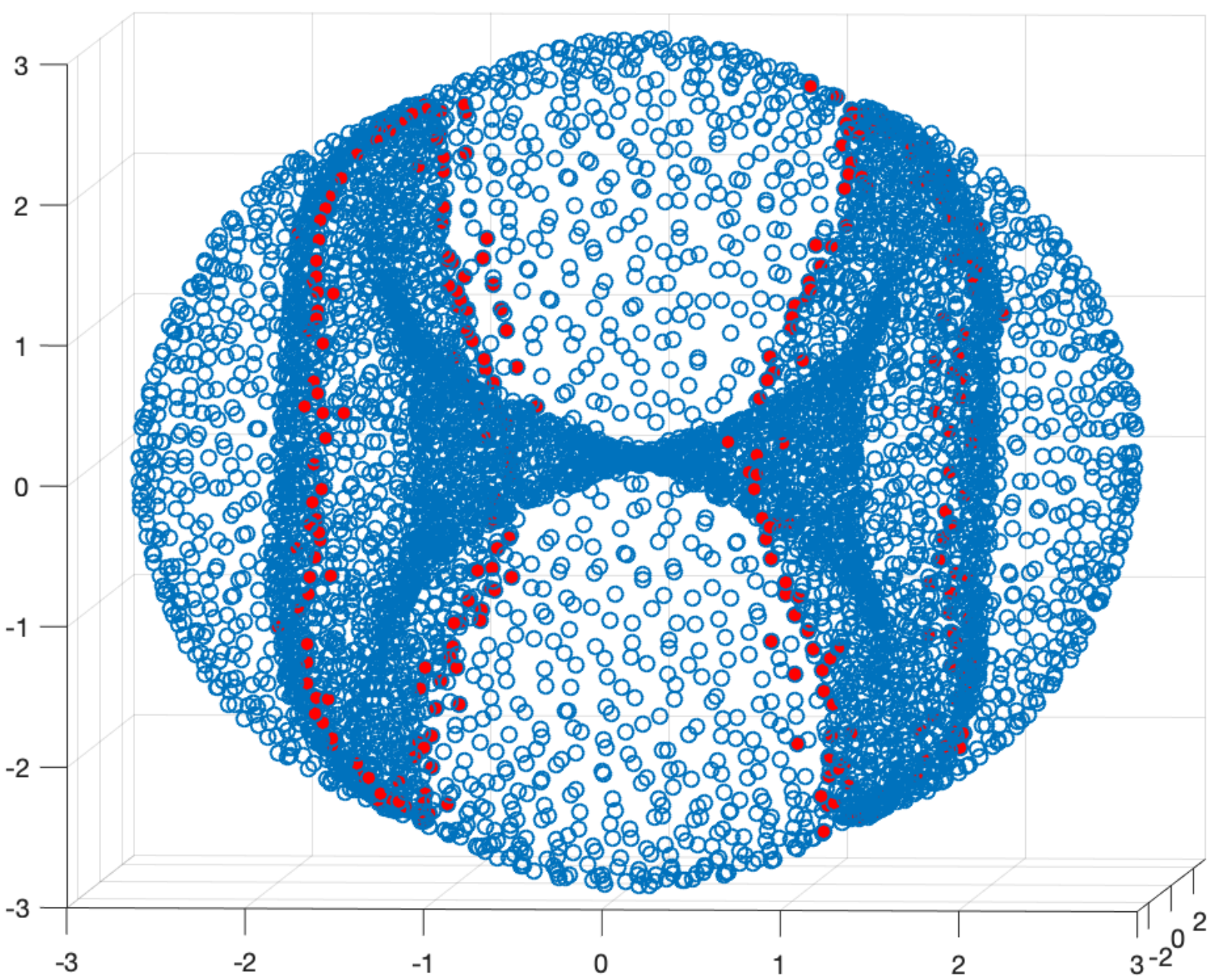}}%
	\subcaptionbox{}{\includegraphics[width=.4\textwidth]{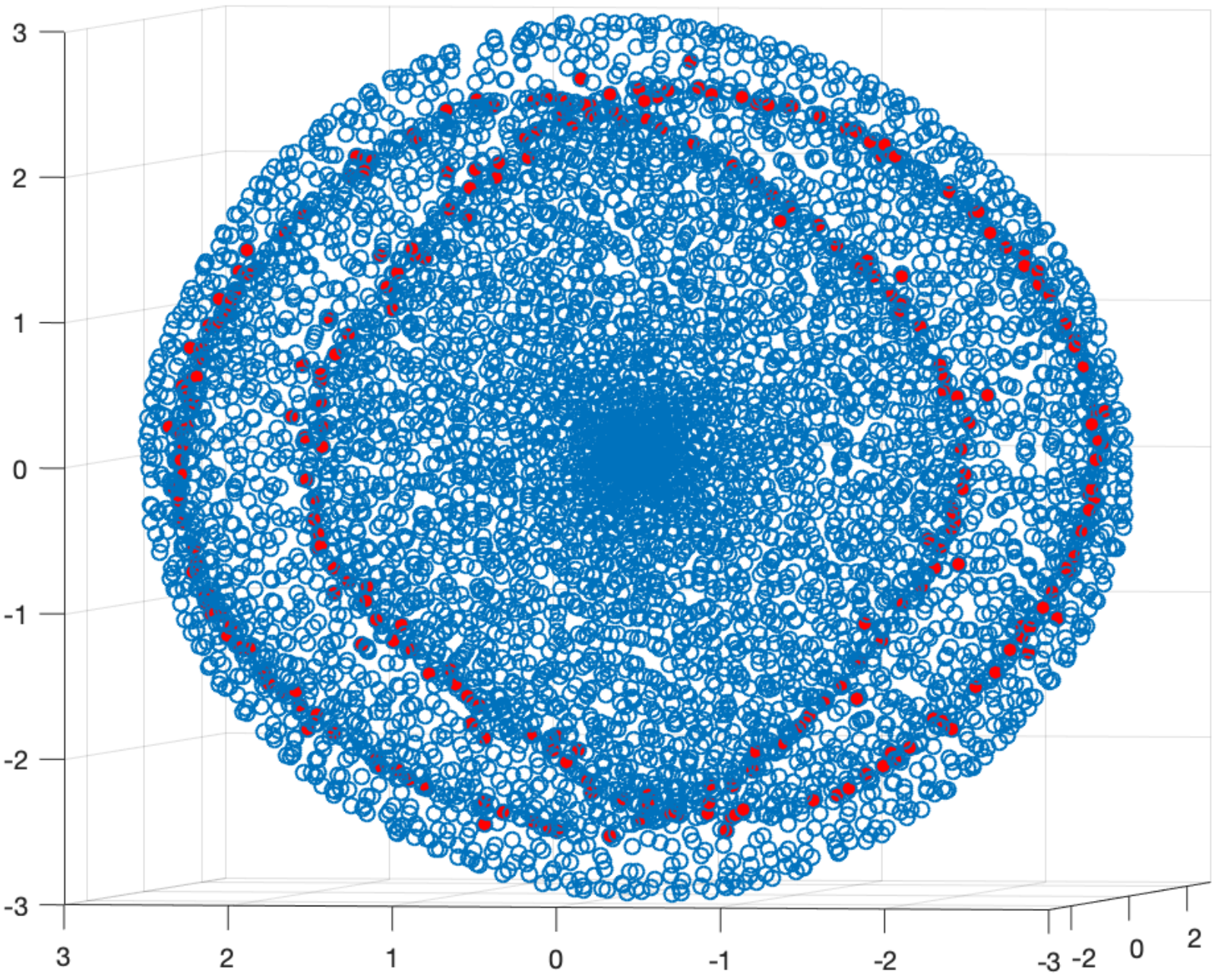}}%
	\caption[Points in the cyclo-octane data set coloured by their local PH properties.]{Two-dimensional {\sc IsoMAP} projection\cite{tenenbaum2000global} of points sampled from the 24-dimensional conformation space of cyclooctane. Points $x$ for which $\dim \Hom^1(A_x) > 1$ have been coloured red, and these clearly appear to cluster near the two embedded circles where the two surfaces intersect.}\label{Fig: Homology Plots Cyclooctane}
\end{figure}

\begin{figure}[ht!]
	\centering
	\subcaptionbox{}{\includegraphics[width=.5\textwidth]{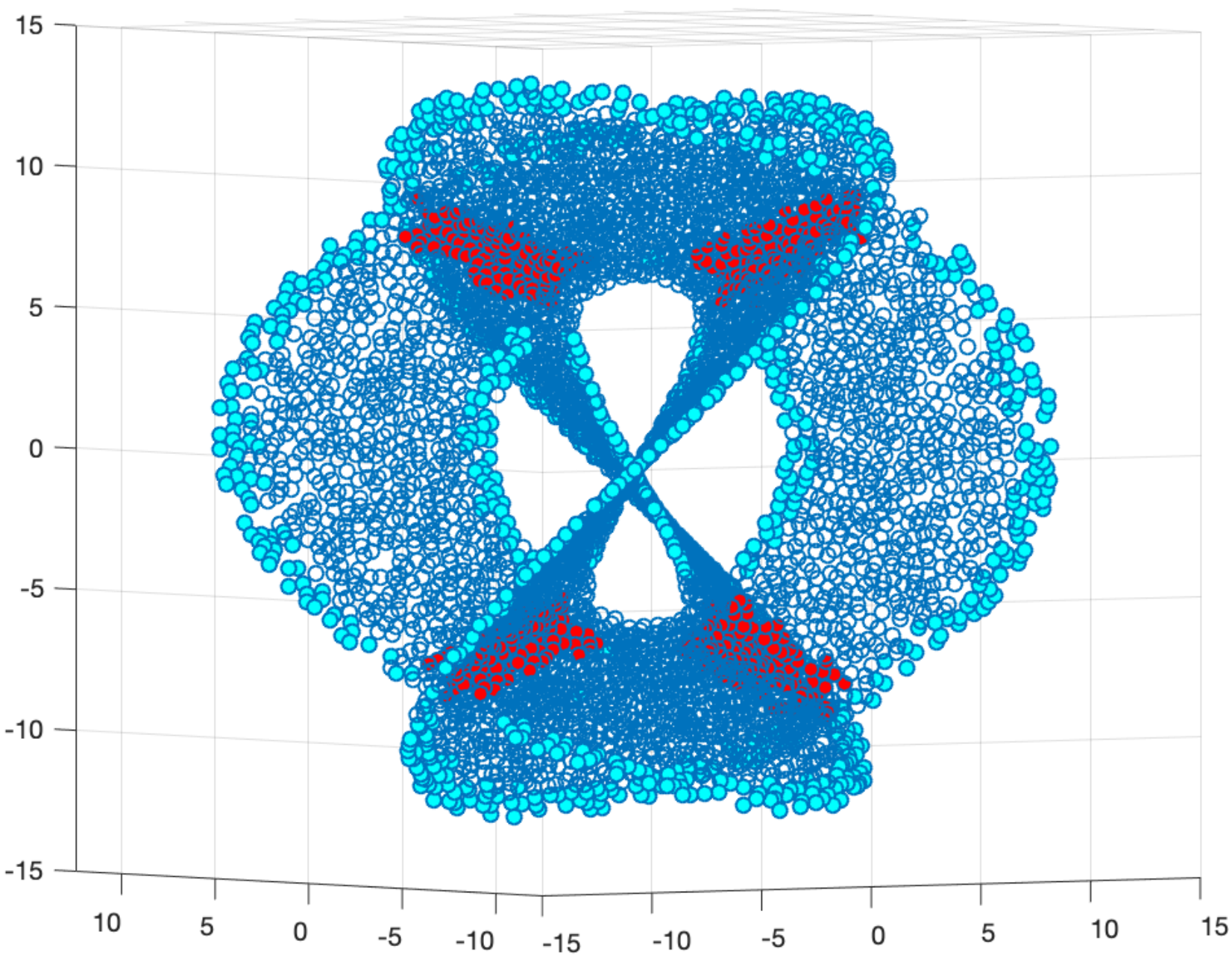}}%
	\subcaptionbox{}{\includegraphics[width=.5\textwidth]{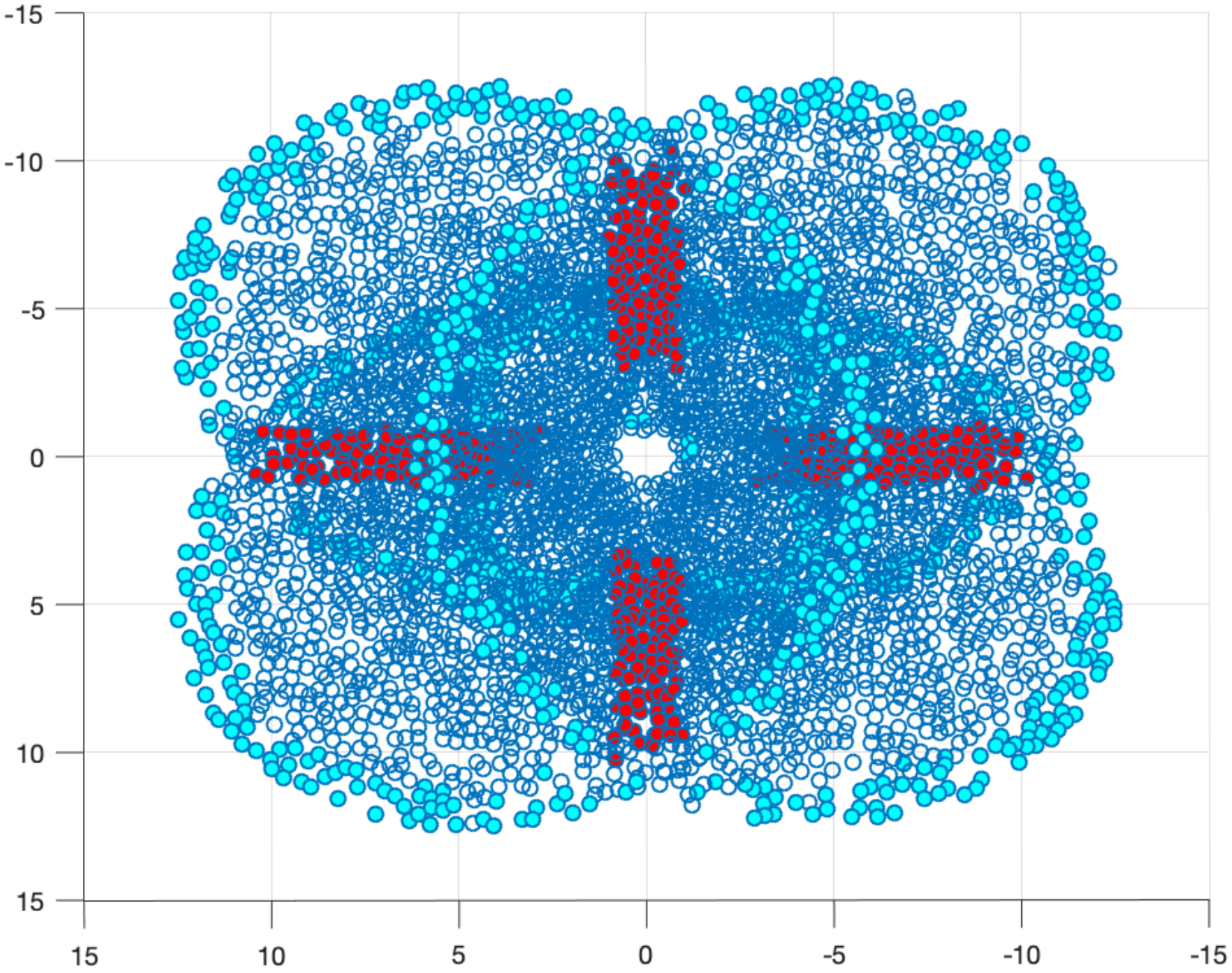}}%
	\caption[Points in the Henneberg surface data set coloured by their local PH properties.]{Two-dimensional projections of points sampled from Henneberg's minimal surface embedded in 3-dimensional space. Points $x$ for which $\dim \Hom^1(A_x) > 1$ are shown in red, and these lie along the four self-intersections. Similarly, points $x$ for which $\dim \Hom^1(A_x) = 0$ have been coloured cyan and appear near the boundary.}\label{Fig: Homology Plots Henneberg}
\end{figure}

%262 words
Local persistent cohomology successfully identifies all the non-manifold regions in two completely different data sets whose underlying spaces are known to admit singularities. The first of these is the {\em conformation space of the cyclo-octane molecule} $\text{C}_{8}\text{H}_{16}$. A single molecule consists of eight carbon atoms arranged in a ring, with each carbon atom being bound to two other carbon atoms and two hydrogen atoms. Under the influence of external chemical and physical forces, cyclo-octane assumes different forms, or {\em conformations}, in 3-dimensional space. The locations of hydrogen atoms are completely determined by those of the carbon atoms, so each conformation may be represented by a point in $\R^{24}$ (i.e., three spatial coordinates for each of the eight carbon atoms). The space of all such conformations forms the union of a Klein bottle and a sphere along two circles\cite{Martin2010,Martin2011}. Having sampled points from this conformation space, we depict (a two-dimensional projection of) the partition of data points by local persistent cohomology in Figure \ref{Fig: Homology Plots Cyclooctane} --- points lying near the two singular circles are indeed separated from all other points. Our second dataset is obtained by uniformly sampling points from the non-orientable {\em Henneberg minimal surface}, which is an immersion of 2-dimensional projective space in standard 3-dimensional space. The results are depicted in Figure \ref{Fig: Homology Plots Henneberg}: again, the points which lie near the four self-intersections are manifestly separated from manifold-like points and boundary points.

% 247 words
Given the enormous quantities of heterogeneous data which are being generated by modern experimental tools, continued reliance on the manifold hypothesis for geometric modelling will become far less tenable in the future. The procedure described here takes the first steps towards relaxing the manifold assumption by enabling us to identify singularities from the local topology of data points. Aside from the data-dependent choice of radius parameters $r$ and $s$ which determine the sizes of annular neighbourhoods $A_x$, the method is entirely unsupervised. Moreover, it enjoys three remarkably convenient properties for our purposes. First, it can be iterated to discover more refined singularities of lower dimension: for instance, had the red points from Figure \ref{Fig: Homology Plots Cyclooctane} formed a singular space of their own (such as a figure-eight rather than disjoint circles), we could have repeated our cohomological clustering operation on the subset of red points to separate out points lying near the lower singularities. Second, the local cohomology computations which form the backbone of this procedure are easily distributed across a host of processors: the persistent homology of annular neighbourhoods $A_x$ and $A_y$ for distinct points $x$ and $y$ in a dataset can -- and should -- be computed in parallel. And third, since persistent cohomology is stable with respect to bounded noise\cite{stability}, the clustering produced by this method inherits a degree of robustness to perturbations of the original dataset.

\begin{methods}
Detailed accounts of the first three topics described below may be found in the textbooks of Hatcher\cite{hatcher}, Oudot\cite{oudot} and Kirwan-Woolf\cite{kirwan2006introduction} respectively.

\subsection{The cohomology of simplicial complexes.}
% 410 words	
	A {\em simplicial complex} $K$ is a collection of subsets of a finite set $V$ (usually called the set of vertices) satisfying the following condition: if $\sigma \subset V$ is in $K$ and $\tau \subset \sigma$ then $\tau$ is also in $K$. The {dimension} of a simplex $\sigma$ is one less than its cardinality, and the set of all $i$-dimensional simplices in $K$ is denoted $K(i)$. The most familiar simplicial complexes are graphs, where $K(0)$ and $K(1)$ correspond to vertices and edges respectively. For each $i$-dimensional simplex $\sigma$, denote by $1_\sigma:K(i) \to \R$ the {characteristic function} which evaluates to $1$ on $\sigma$ and $0$ on all other simplices. The vector space obtained by treating all such characteristic functions as an orthonormal basis is written $\Chn^i(K)$ and called the space of {\em $i$-cochains}. It is possible to construct a sequence of {\em coboundary operators} $\delta^i:\Chn^i(K) \to \Chn^{i+1}(K)$ with the following matrix representation in our chosen basis: the entry in $1_\sigma$'s column and $1_\tau$'s row equals $\pm 1$ if $\sigma \subset \tau$ and is $0$ otherwise. It is always possible to choose signs of the nonzero entries consistently so that the kernel of $\delta^i$ contains the image of $\delta^{i+1}$, and the $i$-th {\em cohomology} of $K$ is the quotient vector space $\Hom^i(K) = \ker \delta^i / \text{img }\delta^{i-1}$. 
	
	Cohomology is an extremely well-studied\cite{hatcher} descriptor of simplicial complexes and related spaces; it enjoys many wonderful properties, but only two of them are relevant to our purposes here. First, it is a {\em homeomorphism invariant}, meaning that any two different triangulations of the same space $X$ will produce identical cohomologies even though the cochain spaces and coboundary operators might be wildly different. For instance, the cohomology vector spaces of an $n$-sphere depend neither on geometric intricacies (such as its radius or its embedding in Euclidean space) nor on the combinatorics of a particular choice of simplicial decomposition. And second, cohomology is {\em functorial} with respect to the subcomplex relation among simplicial complexes. A subset $L$ of simplices in $K$ is called a subcomplex if it  happens to be a simplicial complex in its own right. Whenever $L$ is a subcomplex of $K$, there are well-defined linear maps $\Hom^i(K) \to \Hom^i(L)$ induced on the associated cohomology vector spaces. 
	
\subsection{The persistent cohomology of data.}
% 295 words
	Given a finite dataset $P$ embedded in Euclidean space $\R^n$ and a scale parameter $t \geq 0$, the {\em Vietoris-Rips} simplicial complex $\VR_t(P)$ contains as its $i$-dimensional simplices all subsets $\{p_0,\ldots,p_i\}$ of $P$ whose pairwise Euclidean distances $\|p_j-p_k\|$ are no larger than $t$. It follows that $\VR_t(P)$ is a subcomplex of $\VR_u(P)$ whenever $t \leq u$. By the functoriality of cohomology, in each dimension $i \geq 0$ we obtain not only a one-parameter family of cohomology vector spaces 
	\[
	V(t) = \Hom^i\big(\VR_t(P)\big),
	\]
	but also a compatible family of induced linear maps $V(u) \to V(t)$ for all pairs of real numbers $t \leq u$. Such collections of vector spaces and linear maps indexed by the positive real numbers are called {\em persistence modules}, and their systematic study -- which forms the theoretical core of topological data analysis -- has been greatly facilitated by three miraculous properties.
	
	The first property is algebraic --- although persistence modules appear to involve an infinite amount of information prima facie, any $V$ arising from the Vietoris-Rips cohomology of a finite dataset $P \subset \R^n$ is completely determined by a finite collection $\Barc(P)$ comprising certain half-open subintervals of $\R$, called the {\em barcode} of $P$. The second property is computational; barcodes can be extracted via elementary matrix algebra, and there are several software packages dedicated to their efficient computation\cite{roadmap}. The third crucial property of persistence modules is geometric, and takes the form of a {\em stability theorem}\cite{stability}. Roughly, this result asserts that if the points of $P$ are perturbed by an amount $\epsilon>0$, then the intervals in $\Barc(P)$ also have their endpoints shifted by no more than $\epsilon$. As a consequence, one can conclude that Vietoris-Rips persistent cohomology barcodes are robust to the presence of bounded noise in the original dataset.

\subsection{Stratified Spaces.}

Singular spaces, such as algebraic varieties and quotients of group actions on manifolds, are often analysed via their {\em stratifications}. Each stratification $Y_\bullet$ of an $n$-dimensional space $Y$ is an ascending sequence of closed subspaces
\[
\varnothing = Y_{-1} \subset Y_0 \subset Y_1 \subset \cdots \subset Y_{n-1} \subset Y_n = Y
\]
where the connected components of successive differences $Y_i - Y_{i-1}$, called the {\em $i$-strata}, are open $i$-dimensional submanifolds of $Y$. Every simplicial complex, for instance, admits a natural stratification whose $i$-strata are precisely the $i$-simplices. It is customary to impose two additional constraints on the strata in order to render the study of stratified spaces tractable. The first requirement, called the {\em frontier axiom}, ensures that the set of all strata is partially ordered by the boundary relation $\sigma \leq \tau$ whenever the closure of $\tau$ intersects $\sigma$ (this mirrors the ordering on simplices given by the containment relation $\sigma \subset \tau$). The second requirement, called {\em equisingularity} or {\em normal triviality}, imposes severe topological constraints on intersections of small neighbourhoods in $Y$ around various points of a single $i$-stratum with the higher strata $Y_j$ for $j \geq i$.

As a consequence of equisingularity, to each $i$-stratum $\sigma$ one can assign a single $(n-i-1)$-dimensional stratified space $L_\bullet$, called the {\em link} of $\sigma$, so that the following property holds. For each point $y$ in $\sigma$ and all choices of small neighbourhoods $U_y \subset Y$ of $y$, the intersection of $U_y$ with higher strata $Y_j$ admits a tangent $\times$ normal decomposition of the form
\[
U_y \cap Y_j = \mathbb{R}^i \times \text{Cone}(L_{j-i-1}),
\]
where $\text{Cone}(L_\bullet)$ is the quotient of $L_\bullet \times [0,1)$ obtained by identifying all pairs of the form $(\ell,0)$ with a single point. When $j=i$, we have $U_y \cap \sigma = \mathbb{R}^i$, thus guaranteeing that $\sigma$ is an $i$-dimensional manifold. And for $j=n$, we have $U_y = \mathbb{R}^i \times \text{Cone}(L)$, so it follows that the homeomorphism type -- and hence, the cohomology -- of the boundary $\partial U_y$ is independent of the choice of $y$ in $\sigma$. This is the key property of stratified spaces which is used in our algorithm to identify singular regions within datasets. In this discrete setting, we have no direct access to $\partial U_y$ for a given data point $y$; however, we are able to approximate its cohomology via the persistent cohomology of all the data points lying within an annular neighbourhood $A_y$ of $y$.

\subsection{Datasets.} The {\bf cyclo-octane dataset}, which was introduced by Martin et al.\cite{Martin2010}, consists of 6040 points in $\R^{24}$ subsampled from a far larger dataset containing over a million cyclo-octane conformations. This data set is publicly available as part of the {\sc javaPlex} software package\cite{Javaplex}. The {\bf Henneberg surface dataset} was kindly provided by Martin et al.\cite{Martin2011}; it consists of $5456$ points sampled from the Henneberg surface using the following parametrisation:
		\begin{align*}
			x &= \frac{2(\beta^2 - 1) \cos(\phi)}{\beta} - \frac{2(\beta^6-1) \cos(3\phi)}{3\beta^3},\\
			y &= -\frac{6\beta^2 (\beta^2-1) \sin(\phi) + 2(\beta^6-1) \sin(3\phi)}{3\beta^3},\\
			z &=  \frac{2(\beta^4+1)\cos(2\phi)}{\beta^2},
\end{align*}
where $\beta \in [0.4,0.6]$ and $\phi \in [0,2\pi]$. In this range of $\beta$-values, the surface does not have triple-intersections. 

\subsection{Algorithm and Implementation.}

Procedure \ref{alg:ClassificationPHbars} discovers intersections of dimension $(k-1)$ from points sampled on $k$-dimensional submanifolds of $\R^n$ for $n > k$. It can be suitably generalised in order to find lower-dimensional singularities\cite{nanda2017local}. The partition produced by Procedure \ref{alg:ClassificationPHbars} decomposes the original dataset $P$ into the $k$-manifold points $P_\text{man}$, the boundary points $P_\text{bnd}$ and the desired intersection points $P_\text{int}$. We have implemented Procedure~\ref{alg:ClassificationPHbars} in {\sc Matlab} for surfaces, i.e., for $k=2$, using the inbuilt function {\tt rangesearch} to compute the annuli $A_y$ and {\sc Ripser}\cite{bauer2017ripser} for persistent cohomology barcodes. The annulus parameters $(r,s)$ equal $(0.4,0.25)$ for the cyclo-octane data and $(2,1.5)$ for the Henneberg surface data. The projections of Figure \ref{Fig: Homology Plots Cyclooctane} were obtained by initialising {\sc IsoMAP}\cite{tenenbaum2000global} with $5$-nearest neighbours. 

\begin{algorithm}[h!]
	\caption[Geometric Anomaly Detection]{\em Geometric Anomaly Detection}
	\label{alg:ClassificationPHbars}
	\begin{algorithmic} 
		\REQUIRE Finite point set $P \subset \R^n$, real parameters $0 < r < s$.
		\ENSURE A partition of $P$ into three subsets, $P_\text{man}$, $P_\text{bnd}$ and $P_\text{int}$ 
		\STATE Initialise $P_\text{man}$, $P_\text{bnd}$ and $P_\text{int}$ to $\varnothing$
		\FORALL{$y \in P$}  
		\STATE Find ${A}_{y}\subset P$ containing all $x$ in $P$ which satisfy $r \leq \|x-y\| \leq s$
		\STATE Compute $\Barc_{k-1}(A_y)$, the $(k-1)$-dim Vietoris-Rips barcode of $A_y$ 
		\STATE Calculate $N_y$, the number of intervals in $\Barc_{k-1}(A_y)$ of length $> (s - r)$ 
		\IF{$N_y = 0$}
		\STATE Add $y$ to $P_\text{bnd}$
		\ELSIF{$N_y = 1$}
		\STATE Add $y$ to $P_\text{man}$
		\ELSE{}
		\STATE Add $y$ to $P_\text{int}$
		\ENDIF
		\ENDFOR
	\end{algorithmic}
\end{algorithm}

\end{methods}

\bibliographystyle{naturemag}
\bibliography{lochom}

\begin{thebibliography}{10}
\expandafter\ifx\csname url\endcsname\relax
  \def\url#1{\texttt{#1}}\fi
\expandafter\ifx\csname urlprefix\endcsname\relax\def\urlprefix{URL }\fi
\providecommand{\bibinfo}[2]{#2}
\providecommand{\eprint}[2][]{\url{#2}}

\bibitem{MR3522608}
\bibinfo{author}{Fefferman, C.}, \bibinfo{author}{Mitter, S.} \&
  \bibinfo{author}{Narayanan, H.}
\newblock \bibinfo{title}{Testing the manifold hypothesis}.
\newblock \emph{\bibinfo{journal}{J. Amer. Math. Soc.}}
  \textbf{\bibinfo{volume}{29}}, \bibinfo{pages}{983--1049}
  (\bibinfo{year}{2016}).
\newblock \urlprefix\url{https://doi.org/10.1090/jams/852}.

\bibitem{whatispca}
\bibinfo{author}{Ringner, M.}
\newblock \bibinfo{title}{What is principal component analysis?}
\newblock \emph{\bibinfo{journal}{Nature biotechnology}}
  \textbf{\bibinfo{volume}{26}}, \bibinfo{pages}{303--4}
  (\bibinfo{year}{2008}).

\bibitem{Seung2268}
\bibinfo{author}{Seung, H.~S.} \& \bibinfo{author}{Lee, D.~D.}
\newblock \bibinfo{title}{The manifold ways of perception}.
\newblock \emph{\bibinfo{journal}{Science}} \textbf{\bibinfo{volume}{290}},
  \bibinfo{pages}{2268--2269} (\bibinfo{year}{2000}).
\newblock \urlprefix\url{https://science.sciencemag.org/content/290/5500/2268}.

\bibitem{5714408}
\bibinfo{author}{{Vidal}, R.}
\newblock \bibinfo{title}{Subspace clustering}.
\newblock \emph{\bibinfo{journal}{IEEE Signal Processing Magazine}}
  \textbf{\bibinfo{volume}{28}}, \bibinfo{pages}{52--68}
  (\bibinfo{year}{2011}).

\bibitem{DBLP:conf/iclr/CheLJBL17}
\bibinfo{author}{Che, T.}, \bibinfo{author}{Li, Y.}, \bibinfo{author}{Jacob,
  A.~P.}, \bibinfo{author}{Bengio, Y.} \& \bibinfo{author}{Li, W.}
\newblock \bibinfo{title}{Mode regularized generative adversarial networks}.
\newblock In \emph{\bibinfo{booktitle}{5th International Conference on Learning
  Representations, {ICLR} 2017, Toulon, France, April 24-26, 2017, Conference
  Track Proceedings}} (\bibinfo{year}{2017}).
\newblock \urlprefix\url{https://openreview.net/forum?id=HJKkY35le}.

\bibitem{NIPS2017_6975}
\bibinfo{author}{Sabour, S.}, \bibinfo{author}{Frosst, N.} \&
  \bibinfo{author}{Hinton, G.~E.}
\newblock \bibinfo{title}{Dynamic routing between capsules}.
\newblock In \bibinfo{editor}{Guyon, I.} \emph{et~al.} (eds.)
  \emph{\bibinfo{booktitle}{Advances in Neural Information Processing Systems
  30}}, \bibinfo{pages}{3856--3866} (\bibinfo{publisher}{Curran Associates,
  Inc.}, \bibinfo{year}{2017}).
\newblock
  \urlprefix\url{http://papers.nips.cc/paper/6975-dynamic-routing-between-capsules.pdf}.

\bibitem{nanda2017local}
\bibinfo{author}{Nanda, V.}
\newblock \bibinfo{title}{Local cohomology and stratification}.
\newblock \emph{\bibinfo{journal}{Foundations of Computational Mathematics}}
  (\bibinfo{year}{2019}).

\bibitem{kirwan2006introduction}
\bibinfo{author}{Kirwan, F.} \& \bibinfo{author}{Woolf, J.}
\newblock \emph{\bibinfo{title}{An introduction to intersection homology
  theory}} (\bibinfo{publisher}{Chapman and Hall/CRC}, \bibinfo{year}{2006}).

\bibitem{IH2}
\bibinfo{author}{Goresky, M.} \& \bibinfo{author}{MacPherson, R.}
\newblock \bibinfo{title}{Intersection homology {II}}.
\newblock \emph{\bibinfo{journal}{Inventiones Mathematicae}}
  \textbf{\bibinfo{volume}{71}}, \bibinfo{pages}{77--129}
  (\bibinfo{year}{1983}).

\bibitem{fulton}
\bibinfo{author}{Fulton, W.}
\newblock \emph{\bibinfo{title}{Intersection Theory}}
  (\bibinfo{publisher}{Springer-Verlag}, \bibinfo{year}{1998}).

\bibitem{mischaikow:nanda:13}
\bibinfo{author}{Mischaikow, K.} \& \bibinfo{author}{Nanda, V.}
\newblock \bibinfo{title}{Morse theory for filtrations and efficient
  computation of persistent homology}.
\newblock \emph{\bibinfo{journal}{Discrete and Computational Geometry}}
  \textbf{\bibinfo{volume}{50}}, \bibinfo{pages}{330--353}
  (\bibinfo{year}{2013}).

\bibitem{henselman:ghrist:16}
\bibinfo{author}{Henselman, G.} \& \bibinfo{author}{Ghrist, R.}
\newblock \bibinfo{title}{Matroid filtrations and computational persistent
  homology}.
\newblock \emph{\bibinfo{journal}{arXiv:1606.00199 [math.AT]}}
  (\bibinfo{year}{2016}).

\bibitem{markov}
\bibinfo{author}{Markov, A.~A.}
\newblock \bibinfo{title}{O konstruktivnykh funkciyakh}.
\newblock \emph{\bibinfo{journal}{Trudy Mat. Instituta im. Steklova}}
  \textbf{\bibinfo{volume}{52}}, \bibinfo{pages}{315--348}
  (\bibinfo{year}{1958}).

\bibitem{Martin2011}
\bibinfo{author}{Martin, S.} \& \bibinfo{author}{Watson, J.-P.}
\newblock \bibinfo{title}{Non-manifold surface reconstruction from
  high-dimensional point cloud data}.
\newblock \emph{\bibinfo{journal}{Computational Geometry}}
  \textbf{\bibinfo{volume}{44}}, \bibinfo{pages}{427--441}
  (\bibinfo{year}{2011}).

\bibitem{tenenbaum2000global}
\bibinfo{author}{Tenenbaum, J.~B.}, \bibinfo{author}{De~Silva, V.} \&
  \bibinfo{author}{Langford, J.~C.}
\newblock \bibinfo{title}{A global geometric framework for nonlinear
  dimensionality reduction}.
\newblock \emph{\bibinfo{journal}{science}} \textbf{\bibinfo{volume}{290}},
  \bibinfo{pages}{2319--2323} (\bibinfo{year}{2000}).

\bibitem{Martin2010}
\bibinfo{author}{Martin, S.}, \bibinfo{author}{Thompson, A.},
  \bibinfo{author}{Coutsias, E.~A.} \& \bibinfo{author}{Watson, J.-P.}
\newblock \bibinfo{title}{Topology of cyclo-octane energy landscape}.
\newblock \emph{\bibinfo{journal}{The Journal of Chemical Physics}}
  \textbf{\bibinfo{volume}{132}}, \bibinfo{pages}{234115}
  (\bibinfo{year}{2010}).

\bibitem{stability}
\bibinfo{author}{Cohen-Steiner, D.}, \bibinfo{author}{Edelsbrunner, H.} \&
  \bibinfo{author}{Harer, J.}
\newblock \bibinfo{title}{Stability of persistence diagrams}.
\newblock \emph{\bibinfo{journal}{Discrete and Computational Geometry}}
  \textbf{\bibinfo{volume}{37}}, \bibinfo{pages}{107--120}
  (\bibinfo{year}{2007}).

\bibitem{hatcher}
\bibinfo{author}{Hatcher, A.}
\newblock \emph{\bibinfo{title}{Algebriac Topology}}
  (\bibinfo{publisher}{Cambridge University Press}, \bibinfo{year}{2002}).

\bibitem{oudot}
\bibinfo{author}{Oudot, S.}
\newblock \emph{\bibinfo{title}{Persistence Theory: from Quiver Representations
  to Data Analysis}} (\bibinfo{publisher}{American Mathematical Society},
  \bibinfo{year}{2015}).

\bibitem{roadmap}
\bibinfo{author}{Otter, N.}, \bibinfo{author}{Porter, M.~A.},
  \bibinfo{author}{Tillmann, U.}, \bibinfo{author}{Grindrod, P.} \&
  \bibinfo{author}{Harrington, H.~A.}
\newblock \bibinfo{title}{A roadmap for the computation of persistent
  homology}.
\newblock \emph{\bibinfo{journal}{EPJ Data Science}}
  \textbf{\bibinfo{volume}{6}} (\bibinfo{year}{2017}).
\newblock \urlprefix\url{https://doi.org/10.1140/epjds/s13688-017-0109-5}.

\bibitem{Javaplex}
\bibinfo{author}{Tausz, A.}, \bibinfo{author}{Vejdemo-Johansson, M.} \&
  \bibinfo{author}{Adams, H.}
\newblock \bibinfo{title}{Java{P}lex: {A} research software package for
  persistent (co)homology}.
\newblock In \bibinfo{editor}{Hong, H.} \& \bibinfo{editor}{Yap, C.} (eds.)
  \emph{\bibinfo{booktitle}{Proceedings of ICMS 2014}}, Lecture Notes in
  Computer Science 8592, \bibinfo{pages}{129--136} (\bibinfo{year}{2014}).
\newblock \bibinfo{note}{Software available at
  \url{http://appliedtopology.github.io/javaplex/}}.

\bibitem{bauer2017ripser}
\bibinfo{author}{Bauer, U.}
\newblock \bibinfo{title}{Ripser: a lean {\sc c++} code for the computation of
  vietoris--rips persistence barcodes}.
\newblock \bibinfo{howpublished}{Software available at \url{https://github.
  com/Ripser/ripser}} (\bibinfo{year}{software retrieved in 2017}).

\end{thebibliography}

\begin{addendum}
 \item We thank Barbara Mahler for performing the isomap projection of the cyclo-octane data to $\R^3$. BJS thanks the EPSRC and MRC (EP/G037280/1) and F. Hoffmann-La Roche AG for funding her doctoral studies. 
HAH acknowledges funding from a
Royal Society University Research Fellowship. VN’s work was supported by The Alan Turing Institute under the EPSRC grant
	number EP/N510129/1.
 \item[Correspondence] Correspondence and requests for materials
should be addressed to Bernadette J Stolz~(email: stolz@maths.ox.ac.uk).
\end{addendum}

\end{document}